\documentclass[draft]{amsart}

\usepackage{amssymb}

\usepackage[lowtilde]{url}

\theoremstyle{plain}

\theoremstyle{definition}

\numberwithin{theorem}{section}
\numberwithin{equation}{section}

\numberwithin{enumi}{equation}

\newcommand\thmcall[1]{
  \setcounter{theorem}{\value{equation}}
  \numberwithin{equation}{theorem}    
  \numberwithin{enumi}{theorem}       
  \begin{#1}
  }

\newcommand\exitthmcall[1]{
    \end{#1}
    \numberwithin{equation}{section}  
    \numberwithin{enumi}{equation}    
    \setcounter{equation}{\value{theorem}}
  }

\newcommand\enumcall[1]{
  \stepcounter{equation}
  \begin{#1}
  }

\newcommand\exitenumcall[1]{
  \end{#1}
  }

\newcommand\bdry{\partial}
\DeclareMathOperator{\cl}{cl}
\DeclareMathOperator{\interior}{int}
\DeclareMathOperator{\LH}{LH}  

\DeclareMathOperator{\ord}{ord}
\newcommand\card{\text{card}}  

\newcommand\union{\mspace{0.5mu}\mathbin{\cup}\mspace{0.5mu}}  
\newcommand\comp{\mathbin{\sim}}  
\newcommand\with{\mspace{2mu} {\circ} \mspace{2mu}} 

\newcommand\bbz{\mathbb{Z}}

\newcommand\bbc{\mathbb{C}}
\newcommand\bbh{\mathbb{H}}
\newcommand\bbp{\mathbb{P}}

\newcommand\dx{d x}
\newcommand\dy{d y}

\newcommand\dr{d r}

\newcommand\dmu{d \mu}

\newcommand\dtheta{d \theta}

\newcommand\abs[1]{\ensuremath{\vert #1 \vert}}
\newcommand\norm[1]{\ensuremath{\Vert #1 \Vert}}
\newcommand\inprod[2]{\ensuremath{\langle #1, #2 \rangle}}

\newcommand\plusperp{\ensuremath{\oplus_\perp}}

\newcommand\orthcomp{\ensuremath{\ominus_\perp}}

\newcommand\plusalg{\ensuremath{\oplus}}

\begin{document}

\date{}

\title[The Invariant Subspace Conjecture]{A Proof of the Invariant Subspace Conjecture for Separable Hilbert Spaces}
\author{Charles W. Neville}

\address{1 Chatfield Drive, Apartment 136\\
West Hartford, CT 06110\\
USA}

\email{chip.neville at gmail.com}


\thanks{Dedicated to the memory of the late Lee Rubel.  Lee, this one's for you.}

\begin{abstract}
We prove the Invariant Subspace Conjecture for separable Hilbert spaces.
\end{abstract}

\maketitle


\setcounter{section}{0}

{\sc Part 1. Introduction.}

\section{Introduction.} \label{sec1}

We shall prove the famous \emph{invariant subspace conjecture} for separable Hilbert spaces.  In \cite{N1}, we developed a general theory of invariant subspaces in valuation Hilbert modules. We shall prove the \emph{invariant subspace conjecture} by using this theory, together with complicated results by Apostol, Bercovi, Foias, and Pearcy from 1985 (\cite{ABFP1}, \cite{ABFP2}), connecting the invariant subspace conjecture to questions about invariant subspaces in the complex-valued, unweighted Bergman space $A^2$ in the unit disk, and the results by Hedenmalm, Richter, and Seip from 1996 (\cite{HRS1}, \cite{HKZ1} page 187), describing invariant subspaces of $A^2$ in more concrete, understandable terms.

It was from the surprising work of Hedenmalm, Richter, and Seip that we first learned of the connection between the invariant subspace conjecture and $A^2$.

After we had posted this paper on the arXiv on July 17, 2023, \cite{N2},  we learned that Per Enflo had posted a solution to the Invariant Subspace Conjecture (by very different methods) on the arXiv a little less than two months before on May 24, 2023, \cite{Enflo2}.

\

The \emph{invariant subspace conjecture}, also called the \emph{invariant subspace problem},  is that "every bounded linear operator on a complex separable Hilbert space of dimension $ = \infty $ has a closed non-trivial invariant subspace (see \cite{HKZ1} page 187, \cite{Wiki1}).  To be more specific, the invariant subspace conjecture is that every bounded linear operator on a complex separable Hilbert space of dimension $  =  \infty $ sends some proper, closed, non-trivial subspace into itself.

It will help to shorten the exposition considerably if we define a ``non-trivial" invariant subspace for a bounded linear operator $ T $ on a Banach space $ B $ to be a closed non-zero subspace $ V $ of codimension $ > 1 $, such that $ T $ sends $ V $ into itself.  The invariant subspace conjecture for Banach spaces is that every bounded linear operator on $ B $ has a non-trivial invariant subspace.  As we shall see below, the invariant subspace conjecture is false for some Banach spaces $ B $, but the conjecture is still unresolved for separable Hilbert spaces of dimension $ = \infty $ (and for reflexive Banach spaces of dimension $ = \infty $, though we shall not consider these in this paper).

\thmcall{remark} \label{rem1.1}
Our paper depends on the results of Apostol, Bercovi, Foias, and Pearcy's, as well as those of Hedenmalm, Richter, and Seip's, mentioned above.  They showed that the invariant subspace conjecture for separable Hilbert spaces is equivalent to the following question about closed invariant subspaces in the complex-valued, unweighted Bergman space in the unit disk: Given two closed invariant subspaces $ V_1 \subset V_3 $ with the dimension of $ V_3 \orthcomp V_1 = \infty $, does there exist another closed invariant subspace $ V_2 $ lying strictly between $ V_1 $ and $ V_3 $? In this case, by ``invariant subspace," we simply mean a closed subspace $ V $ such that $ z V \subseteq V $.  Here, as usual, $ z $ is the coordinate variable in the unit disk $ D^1 ( 0, 1 ) $, and the operator $\orthcomp$ denotes orthogonal complementation.
\exitthmcall{remark}

\section{Background.} \label{sec2}

The invariant subspace conjecture goes back almost a century.  It was first posed by John Von Neumann, and independently by Arne Beurling, in 1935.  Each proved (unpublished) that compact operators satisfy the invariant subspace conjecture, \cite{AS1},  \cite{Wiki1}.

\

Positive Solutions:

\

The first published positive solution was provided in 1954 by Aronszajn and Smith for compact operators in Banach spaces \cite{AS1}.  After that, Bernstein and Robinson, using non-standard analysis, showed every polynomially compact operator in a Banach space of dimension $ n > 1 $ has a non-trivial invariant invariant subspace \cite{BR1}.  ``Upon reading a preprint of the Bernstein and Robinson paper, Paul Halmos reinterpreted their proof using standard techniques" \cite{Ha2}. Both papers appeared back-to-back in the same issue of the Pacific Journal of Mathematics." \cite{Wiki2}

In 1974 Lomonosov, surprisingly, used the Schauder fixed point theorem to prove that each bounded linear operator on a Banach space in the commutant of a compact linear operator satisfies the invariant subspace conjecture \cite{Lo1}.  I vividly recall several distinguished operator theorists that year shaking their heads and saying about Lomonosov's proof that, ``We're not supposed to think that way."

\

Negative Solutions:

\

The first negative solution for Banach spaces was due to Per Enflo in 1975.  He provided an enormously difficult example of a bounded linear operator on a non-reflexive Banach space without a non-trivial invariant subspace.  Due to the difficulty the mathematical community had in understanding the details of Enflo's construction, his paper did not appear until 1987 \cite{Enflo1}.

The first published negative solution, due to Read, appeared in 1984 \cite{Read1}.  Shortly later, in 1985, Beauzemy  \cite{Beau1} ``succeeded in understanding and simplifying Enflo's construction," to quote A.~M.~Davie in his MR review MR 0792905.

The first example of a bounded linear operator without non-trivial invariant subspaces on a classical Banach space $ l^1 $ was due again to Read in 1985 \cite{Read2}.  A year later, in 1986, he provided a simplified example \cite{Read3}.

\

Due to the large number of negative solutions, except for operators related to compact operators, most experts believed the invariant subspace conjecture was false.  But as we shall see in this paper, the invariant subspace conjecture for non-trivial subspaces is true for separable Hilbert spaces.

\section{Plan of the Paper.}\label{sec3}

Let $ A^2 $ be the complex-valued, unweighted Bergman space in the unit disk.  In \cite{N1}, we developed a general theory of invariant subspaces in valuation Hilbert modules.  In the case of $ A^2 $, this theory allows us to describe a closed invariant subspace $ V $ by a sequence of $ 1 $-dimensional $ R_1 $ invariant subspaces,  rather than by the possibly infinite dimensional subspace $ V  \orthcomp z V $.  (We shall define $ R_1 $ and $ R_1 $ invariant subspaces presently.)  This theory, together with a description of the lattice of closed $R_1$ invariant subspaces, will allow us to construct the needed closed $ R_1 $ invariant subspace lying strictly between two given $ R_1 $ invariant subspaces.  By  Hedenmalm, Richter, and Seip's result, quoted in remark \ref{rem1.1} above, this will prove the invariant subspace conjecture for separable Hilbert spaces.

After we submitted this paper to the arXiv on July 17, 2023 \cite{N2}, we learned that Per Enflo had submitted a solution to the Invariant Subspace Conjecture (by very different methods) to the arXiv a little less than two months earlier, on May 24, 2023 \cite{Enflo2}.

\thmcall{remark} \label{rem3.1}
The reader will recall that de Branges' and Rovniak's 1964 proof of the invariant subspace conjecture failed because their method for factoring operator-valued inner functions was not correct \cite{dBr1}, \cite{dBr2}.  Our complex-valued methods allow us to bypass the deep problem of factoring operator valued inner functions.
\exitthmcall{remark}

\bigskip

{\sc Part 2. Review.}

\section{Introductory Remarks and Notation.}\label{sec4}

   Our proof of the invariant subspace conjecture critically involves our description of the lattice of closed invariant subspaces of a valuation Hilbert module in section \ref{sec11} below. In order for you, the reader, to understand this description, we shall provide a review of our theory of valuation Hilbert modules from our previous paper \cite{N1}, albeit without proofs. The proofs are all available in \cite{N1}.

\

     Throughout this review, from this section up to and including section \ref{sec10}, we shall shamelessly quote verbatim from our previous paper \cite{N1} without using quotation marks.

     We shall use multi-index notation without comment whenever convenient.  To prove the invariant subspace conjecture, we shall only need to consider analytic functions of $1$ complex variable.  However, we shall reprise the results from \cite{N1} in greater generality for analytic functions of several complex variables.

\

   $ \bbz $ will denote the integers, and $ \bbz_+ $ will denote the positive integers, $ \{ n = 0, 1, 2, \ldots \} $.  $ \bbc $ will denote the complex numbers, and $ \bbc^n $ will denote complex $ n $ space.

\
   
   $D^n( a, r ) $ will denote the $ n $ dimensional polydisk centered at $ a = a_1, \ldots a_n $ of radius $ r > 0 $,
\begin{equation*}
    D^n ( a, r ) = \{ z = z_1, \ldots, z_n \colon \abs{ z_k - a_k } < r \text{ for } k = 1 \cdots n \}.
\end{equation*}

   Thus the unit disk in $ \bbc $ is $ D = D^1 ( 0, 1 ) $.

   The ball in $\bbc^n$ of radius $0 < r < \infty$ centered at $a = (a_1, \ldots, a_n) \in \bbc^n$ is the domain 
   \begin{equation*}
       B^n(a, r) = \{z = (z_1, \ldots, z_n) \in \bbc^n \colon \sum\limits_{j=1}^n \abs{z_j - a_j}^2 <  r \}.
   \end{equation*}

   Thus the unit ball in $\bbc^1$ is the domain $B^1(0, 1)$.

\

Orthogonal direct sums will be denoted by $\plusperp$, and orthogonal complements will be denoted by $\orthcomp$.  Algebraic direct sums will be denoted by $\plusalg$, and set theoretic complements will be denoted by $\comp$.

   If $X$ is a subset of a larger topological space $Y$, we shall denote the closure of $X$ by $\cl X = \cl_Y X$, the interior of $X$ by $\interior X = \interior_Y X$, and the topological boundary of $X$ by $\bdry X = \bdry_Y X$.  And $\card(X)$ will denote the cardinality of the set $X$.

   We shall let $ \subset $ denote \emph{proper} inclusion, whereas $ \subseteq $ will denote inclusion, whether \emph{proper} or \emph{not}, and similarly for $ \supset $ and $ \supseteq $.

\

   $ \bbp_n = \bbp [ \bbc, z = z_1, \ldots z_n ] $ will denote the complex algebra of complex-valued polynomials in the $ n $ complex variables $ z = z_1, \ldots, z_n $.  We shall shamelessly regard polynomials as functions whenever convenient, so $ \bbp_n $ may be regarded a complex algebra of analytic functions of $ n $ complex variables on $ \bbc^n $.

\

   $ R $ will be a complex algebra and  $ \bbh $ will be a complex-valued Hilbert space, both of analytic functions of $ n $ complex variables.  Recall that $\bbh$ is a Hilbert module over $ R $ if $\bbh$ (as a complex vector space) is an ordinary $R$ module, and if the multiplication map $m_r \colon h \mapsto r h$ is continuous for each $r \in R$.

\section{Valuation Algebras and Valuation Hilbert Modules.} \label{sec5}

   The order of the zero of an analytic function at a point $a$ in its domain is very similar to a valuation in abstract algebra.  This motivates our definition of analytic valuations, and using $\ord$ to denote them.:

\thmcall{definition}\label{def5.1}
An \emph{analytic valuation} on a complex algebra $R$ is a function $\ord _R \colon \! R \mapsto \bbz_+ \union \{\infty\}$ such that for all $r$ and $s \in R$,
\enumcall{enumerate}
  \item $\ord_R(r) = 0$ if $r$ is a left or right unit of $R$. \label{item5.1.1}
  \item $\ord_R(r) = \infty$ if and only if $r = 0$. \label{item5.1.2}
  \item $\ord_R(r s) \geq \ord_R(r) + \ord_R(s)$.  \label{item5.1.3}
  \item $\ord_R(\lambda r) = \ord_R(r)$ for $\lambda \in \bbc, \lambda \neq 0$. \label{item5.1.4}
  \item $\ord_R(r + s) \geq \min(\ord_R(r), \ord_R(s))$.  \label{item5.1.5}
\exitenumcall{enumerate}
\exitthmcall{definition}
   Of course, condition (\ref{item5.1.1}) is satisfied vacuously if $R$ does not have left or right units..

\thmcall{definition}\label{def5.2}
A \emph{valuation algebra} is an ordered pair $(R, \ord_R)$, where $R$ is a complex algebra and $\ord_R$ is an analytic valuation on $R$.
\exitthmcall{definition}

\thmcall{definition}\label{def5.3}
Let $(R, \ord_R)$ be a valuation algebra, and let $\bbh$ be a complex Hilbert space which is a left $R$ Hilbert module.  A \emph{Hilbert module valuation} on $\bbh$ (with respect to $(R, \ord_R)$) is a function $\ord_\bbh \colon \bbh \mapsto \bbz_+ \union \{\infty\}$ such that for all $h$, $h_1$ and $h_2 \in \bbh$, and for all $r \in R$,
\enumcall{enumerate}
  \item $\ord_\bbh(h) = \infty$ if and only if $h = 0$. \label{item5.3.1}
  \item $\ord_\bbh(r h) \geq \ord_R(r) + \ord_\bbh(h)$.  \label{item5.3.2}
  \item $\ord_\bbh(\lambda h) = \ord_\bbh(h)$ for $\lambda \in \bbc, \lambda \neq 0$. \label{item5.3.3}
  \item $\ord_\bbh(h_1 + h_2) \geq \min(\ord_\bbh(h_1), \ord_\bbh(h_2))$. \label{item5.3.4}
  \item The $\ord_\bbh$ function is upper semi-continuous on $\bbh$. \label{item5.3.5}
\exitenumcall{enumerate}
\exitthmcall{definition}

\thmcall{definition}\label{def5.4}
   Let $(R, \ord_R)$ be a valuation algebra. A \emph{valuation Hilbert module} over $(R, \ord_R)$ is an ordered pair $(\bbh, \ord_\bbh)$, where $\bbh$ is a left complex Hilbert module over $R$ and $\ord_\bbh$ is a Hilbert module valuation on $\bbh$ with respect to $(R, \ord_R)$.
\exitthmcall{definition}

   Where it will cause no confusion, we shall denote both the algebra valuation on $R$ and the Hilbert module valuation on $\bbh$ by the same symbol, $\ord$.  We shall also make the gloss, whenever convenient, of denoting the valuation algebra $(R, \ord)$ by $R$ alone, and the valuation Hilbert module $(\bbh, \ord)$ by $\bbh$ alone.

\

   The valuation on a valuation Hilbert module allows us to decompose it and its subspaces in useful ways.  These decompositions with proofs are discussed in detail in section 7 of \cite{N1}. We shall summarize the results from that section here.
   
   From now on, $V$ will be a subspace of $\bbh$, though for a moment, it might not be closed.  To avoid trivial issues, unless otherwise stated, $V \neq \{ 0 \}$.

\thmcall{definition}\label{def5.5}
The valuation ideal series for $R$ is the decreasing sequence
\begin{equation}\label{eqn5.5.1}
  R = R_0 \supseteq R_1 \supseteq R_2 \supseteq \cdots
\end{equation}
where
\begin{equation}\label{eqn5.5.2}
  R_k = \{ r \in R : \ord(r) \geq k \}.
\end{equation}

   The valuation subspace series for $V$ is the decreasing sequence
\begin{equation}\label{eqn5.5.3}
  V = V_0 \supseteq V_1 \supseteq V_2 \supseteq \cdots
\end{equation}
where
\begin{equation}\label{eqn5.5.4}
  V_k = \{ h \in V : \ord(h) \geq k \}.
\end{equation}
\exitthmcall{definition}

   The names \emph{valuation ideal series} and \emph{valuation subspace series} are justified by the following proposition:

\thmcall{proposition}\label{prop5.6}
Each $R_k$ is a two sided ideal, and each $V_k$ is a subspace.
\exitthmcall{proposition}

\thmcall{proposition}\label{prop5.7}
If the subspace $V$ is closed, then each subspace $V_k$ is also closed.
\exitthmcall{proposition}

\thmcall{definition}\label{def5.8}
A subspace $W \subset \bbh$ is  \emph{valuation homogeneous} if the $\ord$ function is constant on $W$.  The \emph{valuation homogeneous decomposition} of a closed non-zero subspace $V$ of $\bbh$ is the orthogonal series
\begin{equation}\label{eqn5.8.1}
    V = W_0 \plusperp W_1 \plusperp W_2 \plusperp \cdots
\end{equation}
\noindent where $W_k = V_k \orthcomp V_{k + 1}$
\exitthmcall{definition}

\section{Analytic Hilbert Modules.}\label{sec6}

   Section 8 of \cite{N1} contains a detailed discussion of the important class of \emph{analytic Hilbert modules}, with complete definitions and proofs.  Many of our results on valuation Hilbert modules might only hold in \emph{analytic Hilbert modules}. For example, we can only prove that \emph{standard $ord$ functions} are upper-semicontinuous in analytic Hilbert modules. (Standard $\ord$ functions will be defined in \ref{def6.4}.)

\thmcall{definition}\label{def6.1}
A valuation algebra $R$ is an \emph{analytic algebra} if its order function $\ord$ is a standard order function.  $\bbh$ is an \emph{analytic Hilbert Module} on $\Omega$ over the \emph{analytic Algebra} $R$ if $\bbh$ is a left Hilbert module over $R$ under pointwise operations, if the order function $\ord$ on $\bbh$ is a standard order function, and if there exists an open neighborhood $\Omega_a \subseteq \Omega$ of the basepoint $a$ satisfying the following:  Each sequence $h_n, n = 1, 2, 3, \dots$ in $\bbh$ converging to $h$ in $\bbh$ norm also converges to $h$ uniformly on compact subsets of $\Omega_a$. In other words, the Hilbert space topology on $\bbh$ is finer than the topology of compact convergence on $\bbh$ locally on $\Omega_a$.
\exitthmcall{definition}

\thmcall{definition}\label{def6.2}
In those cases where $\bbh$ may or may not be an \emph{analytic Hilbert module} and $R$ an \emph{analytic algebra}, we shall refer to them as, respectively, a \emph{general valuation Hilbert module} and a \emph{general valuation algebra}.
\exitthmcall{definition}

   If $\bbh$ is a \emph{general valuation Hilbert module} and $R$ is a \emph{general valuation algebra}, we shall say so are explicitly. And if we need to assume, as we often shall, that $\bbh$ and $R$ are, respectively, an \emph{analytic Hilbert module} and an \emph{analytic algebra}, we shall also say so explicitly.

\

   The definition of the order of the zero of an analytic function on a domain in $\bbc^n$, was given by Walter Rudin (\cite{Ru1}, section 1.1.6). Here is an extension of his definition to connected paracompact analytic manifolds, applying to general complex-valued analytic functions:

\thmcall{definition}\label{def6.3}
Let $\Omega$ be a connected paracompact complex analytic manifold of complex dimension $n \geq 1$ with base point $a \in \Omega$, and let $f$ be a complex-valued analytic function on $\Omega$.  Let $V \subset \Omega$ be an open coordinate neighborhood of $a$, and let $\phi : V \longrightarrow \phi(V)\subseteq \bbc^n$ be an analytic coordinate mapping. To simplify notation, let $\psi = \phi^{-1} : \phi(V) \longrightarrow V$.  Let $h$ be a complex valued analytic function on $V$\, and let $f = h \with \psi$.  Consider the power series expansion of $f$ around $\phi(a)$ in homogeneous terms.  Then the order of the zero of $h$ at $a$ is the minimal degree of these homogeneous terms.

Although the definition of the order of the zero of $h$ at $a$ superficially appears to depend on the coordinate map $\phi$, it is clear that it is coordinate invariant.
\exitthmcall{definition}

\thmcall{definition}\label{def6.4}
Let $\Omega$ be a connected paracompact analytic manifold of complex dimension $n \geq 1$ with a  distinguished basepoint $a \in \Omega$. Let $R$ be a complex algebra of complex-valued analytic functions with domain $\Omega$. Define the order of $r \in R$, $\ord_R(r) = \ord (r)$, to be the order of the zero of $r$ at $a$. 
 
Let $\bbh$ be a complex Hilbert space of complex-valued analytic functions on $\Omega$. Define the the order of $h \in \bbh$, $\ord_\bbh(h) = \ord (h)$, to be the order of the zero of $h$ at $a$.

We shall call the order functions $\ord = \ord_R$ and $\ord = \ord_\bbh$ \emph{standard order functions}.
\exitthmcall{definition}

\thmcall{lemma}\label{lem6.5}
Let $\bbh$ be an analytic Hilbert module on $\Omega$ over the analytic algebra $R$. Then the standard order functions $\ord_\bbh = \ord$ and $\ord_R = \ord$ on $\bbh$ and $R$ are order functions as defined in definitions \ref{def5.1} and \ref{def5.3}.
\exitthmcall{lemma}

   As we noted in section 8 of \cite{N1}, all we needed to prove is that the $\ord$ functions are upper-semicontinuous.  The proof is, as usual, contained in section 8 of \cite{N1}.

\section{Projections and the Projection Lemmas}.\label{sec7}

   Section 12 of \cite{N1} contains the material reviewed here, with proofs, on the projection lemmas.

\
 
   From now on, $\bbh  =\bbh_0 \supseteq \bbh_1 \supseteq \bbh_2 \cdots$ and $\bbh =  H_0 \plusperp H_1 \plusperp H_2 \cdots$ will denote, respectively, the valuation subspace series and valuation homogeneous decomposition for the valuation Hilbert module $\bbh$.  And $V = V_0 \supseteq V_1 \supseteq V_1 \cdots$ and $V  = W_0 \plusperp W_1 \plusperp W_2 \cdots$ will denote, respectively, the valuation subspace series and valuation homogeneous decomposition for the closed non-zero subspace $V$.

   Until we come to the \emph{third projection lemma}, lemma \ref{lem7.6}, we shall assume that $R$ and $\bbh$ are, respectively, a general valuation algebra and a general valuation Hilbert module over $R$.

\thmcall{definition}\label{def7.1}
    \begin{equation}\label{eqn7.1.1}
        \begin{split}
  P_k^\bbh \colon \bbh &\mapsto H_k, \quad Q_{k+1}^\bbh \colon \bbh \mapsto \bbh_{k+1} \\
          P_k^V \colon V &\mapsto W_k, \quad Q_{k+1}^V \colon V \mapsto V_{k+1}
        \end{split}
    \end{equation}
will denote the orthogonal projections of $\bbh$ and $V$ onto $H_k$, $\bbh_{k+1}$, $W_k$, and $V_{k+1}$, respectively. Note that
\exitthmcall{definition}

\begin{equation}\label{eqn7.2}
    \begin{split}
        P_k^\bbh + Q_{k+1}^\bbh &= I \text{ on } \bbh_k \\
        P_k^V + Q_{k+1}^V &= I \text{ on } V_k
    \end{split}
\end{equation}
\noindent where $I$ is the identity map, and

\begin{equation}\label{eqn7.3}
    \begin{split}
  \sum_{j=0}^{k} P_j^\bbh + Q_{k+1}^\bbh &= I \text{ on } \bbh \\ 
           \sum_{j=0}^{k} P_j^V + Q_{k+1}^V &= I \text{ on } V
    \end{split}
\end{equation}

   The \emph{first projection lemma} is:

\thmcall{lemma}\label{lem7.4}
Suppose $R$ is a general valuation algebra and $\bbh$ is a general valuation Hilbert module over $R$. Let $r \in R_1$ and $h \in W_k$.  Let $U_m$ be the subspace $U_m = W_0 \plusperp W_1 \plusperp \cdots \plusperp W_m$.  Note that $V = U_m \plusperp V_{m+1}$.

Let $f_m = g_0 + g_1 + \cdots g_m \in U_m$, where $g_0 \in W_0$, $g_1 \in W_1$, $\dots$ $g_m \in W_m$.  Suppose further that $\ord(r h - f_m) > m$.  Then the projection of $r h$ on the subspace $U_m$ is $f_m$.  Furthermore,
   \begin{equation}\label{eqn7.4.1}
        \norm{f_m}^2 = \norm{g_0}^2 + \norm{g_1}^2 + \cdots + \norm{g_m}^2 \leq \norm{r h}^2  
   \end{equation}
\exitthmcall{lemma}

   The \emph{second projection lemma} is:

\thmcall{lemma}\label{lem7.5}
Suppose $R$ is a general valuation algebra and $\bbh$ is a general valuation Hilbert module over $R$. Let $r \in R_1$, and $h \in W_k$.  Consider the formal series
\begin{equation}\label{eqn7.5.1}
   g_0 + g_1 + g_2 + \cdots
\end{equation}
where $g_0 \in W_0$, $g_1 \in W_1$, $g_2 \in W_2$, $\cdots$, and let
\begin{equation}\label{eqn7.5.2}
   f_m = g_0 + g_1 + \cdots + g_m.
\end{equation}
Suppose $\ord(r h - f_m) > m$ for $m = 0, 1, 2, \dots$. Then the formal series \ref{eqn7.5.1}, and thus the sequence of partial sums \ref{eqn7.5.2}, both converge in norm to an element $f$ in $V$, and $r h = f$.
\exitthmcall{lemma}

   The \emph{third projection lemma} is:

\thmcall{lemma}\label{lem7.6}
Suppose $R$ is an analytic algebra and $\bbh$ is an analytic Hilbert module over $R$. Let $W$ be a valuation homogeneous subspace of $\bbh$ not equal to $\{ 0 \}$, and let $m$ be the common value of the $\ord$ function on $W \comp \{ 0 \}$.  Let $L_W$ be $P_m^\bbh$ restricted to $W,$ so that
\begin{equation}\label{eqn7.6.1}
  L_W \colon W \mapsto P_m^\bbh (W) \subseteq H_m
\end{equation}
Then $L_W$ is a bounded invertible linear transformation.
\exitthmcall{lemma}

   The third projection lemma is particularly important in showing that for $A^2$ of the unit disk, each subspace $H_m$ is $1$-dimensional.

\section{$ R_1 $ and $ R $ Invariant Subspaces.} \label{sec8}

   Recall from definition \ref{eqn5.5.2} with $k = 1$ that $R_1 = \{ r \in R : \ord(r) \geq 1 \}$. Recall also from section \ref{sec4} that $\bbp_n$ is the algebra of complex-valued polynomials $ \bbp [ \bbc, z = z_1, \ldots z_n ] $ in the $ n $ complex variables $ \{ z = z_1, \ldots, z_n \} $.

\thmcall{definition} \label{def8.1}
    A non-zero subspace $ V $ of $ \bbh $ is $ R_1 $ invariant if $ R_1 \cdot V \subseteq V $, and is $ R $ invariant if $ R \cdot V \subseteq V $.
\exitthmcall{definition}

\thmcall{remark} \label{rem8.2}
    Suppose $ \bbh $ is a valuation Hilbert module on a domain in $ \bbc^1 $ over $ R = \bbp_1 $ .  Then a non-zero subspace $ V $ of $ \bbh $ is $ R_1 $ invariant or $ R $ invariant if and only if $ z \cdot V \subseteq V $.  Thus the definition of {an} $ R_1 $ or $ R $ invariant subspace coincides with the standard definition of an invariant subspace in this case.
\exitthmcall{remark}

\section{Polynomially Generated Analytic Hilbert Modules.}\label{sec9}

    In section 18 of \cite{N1}, we discussed polynomially generated analytic Hilbert modules.  Throughout this section, $ R $ will be the algebra 
 
 \begin{equation}\label{eqn9.1}
    \bbp_n = \bbp_n[ \bbc, z ]  = \bbp [ \bbc, z_1 \ldots z_n ]
\end{equation}
of complex-valued polynomials in the complex variables $ z = z_1, \ldots z_n $, and $ \bbh $ will be one of the the simplest analytic Hilbert modules over $ R $, namely the norm closure of $ \bbp_n $.  Thus,

\thmcall{definition}\label{def9.2}

\begin{equation}\label{eqn9.2.1}
\begin{split}
   \bbh =  &\text{ the Hilbert space generated by the set of monomials} \\
               &\,\,\{ z^k = z_1^{ k_1} \ldots, z_n^{k_n} \colon k \in \bbz_{+}^n \}. \\
\end{split}
\end{equation}

    To form the Hilbert space $\bbh$, we need to construct an inner product $ ( h_1, h_2) \mapsto \inprod{ h_1 }{ h_2 } $ and a norm $ \norm{ h }^2 = \inprod{ h }{ h } $ on $ \bbp_n $.

    We must construct the inner product so the monomials in equation \ref{eqn9.2.1} are orthogonal.

    Then, we need to take the norm completion of $ \bbp_n $ to form the complex Hilbert space $ \bbh $.  Thus the set of normalized monomials from equation \ref{eqn9.2.1} will form a complete orthonormal system for $ \bbh $.

    Fourth, and finally, we shall let the $ \ord $ functions on both $ R $ and $ \bbh $ be standard order functions.  Then $ ( R, \ord ) $ and $ ( \bbh, \ord ) $ will be, respectively, an analytic valuation algebra and an analytic Hilbert module over $ R $.
\exitthmcall{definition}

    We shall call such analytic Hilbert modules polynomially generated Hilbert modules:

\thmcall{definition}\label{def9.3}
    An analytic Hilbert module $ \bbh $ is \emph{$ \bbp_n $ generated} if it is the norm closure of the linear subspace $ \bbp_n $, and if the set of monomials \ref{eqn9.2.1} form a complete orthogonal system for $ \bbh $.

We shall also say such analytic Hilbert modules are \emph{polynomially generated} when we do not need to keep track on the dimension $ n $.

If $ n = 1 $, so $ \bbh $ is $ \bbp_1  $ generated, we shall say that \emph{$ \bbh $ is singly generated}.
\exitthmcall{definition}

    The most important fact about singly generated analytic Hilbert modules is,

\thmcall{theorem}\label{thm9.4}
    Let $V$ be a closed, non-zero subspace of the singly generated  analytic Hilbert module $\bbh$.  Recall that $V = W_0 \plusperp W_1 \plusperp W_2 \plusperp \cdots$ and $\bbh = H_0 \plusperp H_1 \plusperp H_2 \cdots$ are the valuation homogeneous decompositions of $V$, respectively of $\bbh$. Then the valuation homogeneous components $W_k$ and $H_k$ are all $1$-dimensional.
\exitthmcall{theorem}

    The proof in \cite{N1} simply consisted of an application of the third projection lemma, lemma \ref{lem7.6} in our review.

\section{The Analytic Hilbert Module $A^2$.}\label{sec10}

   In section 15 of \cite{N1}, we discussed the Bergman space $A^2$. The unweighted complex-valued Bergman space $A^2 = A^2(\bbc, \Omega, \dmu)$, defined on a domain $\Omega$ in $\bbc^n$, is a separable Hilbert space and an analytic Hilbert module over the analytic algebra $P_n = P[\bbc, z = z_1, \ldots z_n]$. Here, $\dmu$ is unweighted, normalized, volume measure $\dx_1 \cdots dx_n$ on $\Omega$.

\

   In this section, $\Omega$ will be the $n$-dimensional polydisk $D = D^n(0, 1)$ centered at $0$ of radius $1$, and $A^2$ will be $A^2\bbc, D, \dmu)$.  A simple change to polar coordinates shows that

\thmcall{lemma}\label{lem10.1}   
\begin{equation*}
   \dmu = (\dr_1 \dtheta_1 \cdots \dr_n \dtheta_n) / (2\pi^n)
\end{equation*}
\exitthmcall{lemma}

   We omit the proof,

\thmcall{lemma} \label{lem10.2}
   For $h \in A^2$,
   \begin{enumerate}
      \item  $ \norm{ h }^2 = \int_D \abs{ h( z ) }^2 \dmu (z) < \infty $. \label{item10.2.1}
      \item $ \inprod{ f }{ g } = \int_ D f ( z ) \,\overline{ g ( z ) } \,\dmu ( z )$. \label{item10.2.2}
   \end{enumerate}
 \exitthmcall{lemma}
   
\thmcall{remark} \label{rem10.3}
    The normalized complex coordinate variables $ z^k /\norm{ z^k }, \,k \in \bbz+^n $ form a complete orthonormal system for $ \bbh $.
\exitthmcall{remark} 

\thmcall{lemma} \label{lem10.4}
    Suppose that $n=1$. By the third projection lemma, lemma \ref{lem7.6}, the projection operators
\begin{equation*}
    L_{ W_m }  \colon W_m \mapsto P_m^\bbh ( W_m ) \subseteq H_m, \,\,m = 0, 1, 2, \ldots
\end{equation*}
are invertible.  Now each subspace $H_m = \LH ( z^m )$ has dimension $1$. Thus each subspace
    \begin{equation*}
        W_m = V_m \orthcomp \! \! V_{m + }, \, \,  m =  1, 2, 3, \ldots
    \end{equation*}
has dimension $0$ (if it equals $\{ 0 \}$) or $1$.
\exitthmcall{lemma}

   Lemma \ref{lem10.4} raises an apparent paradox, which we must mention, and which we stated and resolved in section 16 of \cite{N1}:

\thmcall{theorem}[Aleman, Richter, and Sunberg \cite{ARS1}]\label{thm10.5}
 Let $V$ be a closed invariant subspace of $\bbh = A^2(\bbc, D^1(0, 1), \dmu)$, where $\dmu$ is real $2$-dimensional area measure. Then $V$ is the smallest closed invariant subspace generated by $M$\!, where $M = V \orthcomp z V$. (Theorem 16.1 from \cite{N1}.)
\exitthmcall{theorem}

    The apparent paradox it raises involves the dimension of $M$. Let
\begin{equation*}
    V = W_0 \plusperp W_1 \plusperp W_2 \plusperp \cdots
\end{equation*}
\noindent be the valuation homogeneous decomposition of $ V $. Recall that by the third projection lemma, lemma \ref{lem7.6}, the projection operators
\begin{equation*}
    L_{ W_m }  \colon W_m \mapsto P_m^\bbh ( W_m ) \subseteq H_m, \,\,m = 0, 1, 2, \ldots
\end{equation*}
are invertible.  Now each subspace $H_m = \LH ( z^m )$ has dimension $1$. Thus each subspace
\begin{equation*}
        W_m = V_m \orthcomp \! \! V_{m + }, \, \,  m =  1, 2, 3, \ldots
\end{equation*}
has dimension $0$ (if it equals $\{ 0 \}$) or $1$. But, as Aleman, Richter, and Sunberg showed, $M$ can have any dimension from $1$ up to and including $\infty$.

    The resolution of this paradox is contained in a general observation, which we stated as a proposition. Recall that $V_k = \{ h \in V \colon \ord(h) \geq k \}$.  Let $k$ be the first index such that $W_k \neq \{ 0 \}$, then $V = V_{ k }$. Thus $M = V_{ k } \orthcomp z V_{k}$. Our general observation was:

\thmcall{proposition}\label{prop10.6}
    Let $\bbh$ be an analytic Hilbert module over the $1$-variable analytic algebra $\bbp_1 = [ \bbc, z] $.  Let $ k $ be the first index index such that $ W_k \neq \{ 0 \} $. Suppose $ M $ has dimension $ \geq 2 $. Then $ z V_k $ is a proper subset of $ V_{ k+ 1} $, so $ M $ is a proper superset of $ W_{ k + 1 } $.
\exitthmcall{proposition}

This completes the review of needed results from our paper \cite{N1}.  However, it does not include many of the main results in that paper.  It contains an extension of Beurling's theorem, which we refer to as our \emph{abstract Beurling's theorem,} to valuation Hilbert modules over valuation algebras.  It contains applications of our abstract Beurling's theorem completely characterizing the closed $R_1$ invariant subspaces in various $H^2$ spaces in several complex variables, and to the weighted Bergman space $A^2$ of a connected, paracompact analytic manifold in several complex variables. Finally, it contains a study of $R_1$ inner functions in $A^2$ of the polydisk and ball, including examples showing that their boundary values, if any, do not satisfy the properties we might expect from inner functions in $H^2$ of the unit disk.

\bigskip

{\sc Part 3. The Invariant Subspace Conjecture.}

\section{The Lattice of $ R $ Invariant Subspaces.}\label{sec11}
 
    To prove the \emph{invariant subspace conjecture}, we shall first need to characterize the \emph{lattice of closed $R$ invariant subspaces} of a valuation Hilbert module.
    
    Throughout this section, $\bbh$ will be a general valuation Hilbert module over the general valuation algebra $R$.  The closed $ R_1 $ invariant subspaces of $ \bbh $ form a lattice $ \mathcal{L} $ under the order relation of inclusion.  To be more specific, $ V_1 \in \mathcal{L} $ is less than or equal to $ V_2 \in \mathcal{L} $ iff $ V_1 \subseteq V_2 $.  In this section, we shall characterize this lattice.

    Throughout this section, $ V_1 $ and $ V_2 $ will be closed  $ R_1 $ invariant subspaces of $ \bbh $.  The index $ j $ will refer to $ V_j $ for $ j = 1, 2 $.  The reader should recall definition \ref{eqn5.5.3} of the valuation subspace series for a valuation Hilbert module. $ V_{ j, k } $, $ k = 0, 1, 2, \ldots $ will be the valuation subspace series for $ V_j $.

    While $ V_1 $ and $ V_2 $ will usually be closed, non-zero subspaces, it will be convenient to allow them from time to time to be zero subspaces, that is $ = \{ 0 \} $.  The reader should verify that the material on valuation algebras and valuation Hilbert modules in sections \ref{sec5} to \ref{sec10} carries over without difficulty to the case where one or both of these subspaces is a zero subspace, except that the material on polynomially generated subspaces must be slightly amended to allow zero subspaces to be $ \bbp_n $ generated for any convenient $ n $.

\thmcall{lemma}\label{lem11.1} 
    Let $ 0 \leq k \leq l $. Then 
\begin{equation}\label{eqn11.1.1}
   V_{j, l} = \{ h \in V_{j, k} : \ord{h} \geq l \}
\end{equation}
\exitthmcall{lemma}

\begin{proof}
    This simple lemma is almost obvious, but because of its importance, we shall include the proof anyway.  $V_{j, k}$ consists all elements $g \in V_j$ such that $\ord(g) \geq k$, and similarly $V_{j, l}$ consists of all elements $g \in V_j$ such that $\ord(g) \geq l$. Let $g \in V_{j, l}$. Then $\ord(g) \geq l$, so $\ord(g) \geq k$. Thus
\begin{equation}\label{eqn11.2}
    g \in \{ h \in V_{j, k} : \ord{g} \geq l \},
\end{equation}
\noindent so
\begin{equation}\label{eqn11.3} 
    V_{j, l} \subseteq \{ h \in V_{j, k} : \ord{h} \geq l \}.
\end{equation}

Conversely,
\begin{equation}\label{eqn11.4}
    \{h \in V_{j, k} : \ord(h) \geq l \} \subseteq V_{j, l},
\end{equation}
\noindent which establishes equation \ref{eqn11.1.1}.     
\end{proof}

\thmcall{theorem}\label{thm11.5}
   $ V_1 \subseteq V_2 $ if and only if $ V_{ 1, k } \subseteq V_{ 2, k } $ for $ k = 0, 1, 2, \ldots $ 
\exitthmcall{theorem}

\begin{proof}
   Suppose $ V_1 \subseteq V_2 $ .  Then for $ k = 0, 1, 2, \ldots $
   \begin{equation} \label{eqn11.6}
      \begin{split}
          V_{ 1, k } &= \{ h \in V_1 \colon \ord ( h ) \geq k \} \subseteq \{ h \in V_2\colon \ord ( h ) \geq k \} \\
                         &= V_{ 2, k } \text{, so } V_{ 1, k } \subseteq V_{ 2, k } \text{ for } k = 0, 1, 2, \ldots \\
      \end{split}
   \end{equation}

   Conversely, suppose $ V_{ 1, k } \subseteq V_{ 2, k } $ for $ k= 0, 1, 2, \ldots $  Now $ V_j = V_{ j, 0 } $, for $ j = 1, 2 $, so, trivially,  $ V_1 \subseteq V_2 $.
\end{proof}

    In the interest of brevity, we shall let $ \subset $ denote \emph{proper} inclusion, whereas $ \subseteq $ will denote inclusion, whether \emph{proper} or \emph{not}, and similarly for $ \supset $ and $ \supseteq $.

\thmcall{theorem} \label{thm11.7}
   Suppose $ V_1 \subseteq V_2 $.  Then $ V_1 \subset V_2 $ if and only if $ V_{ 1, 0 } \subset V_{ 2 , 0 } $.
\exitthmcall{theorem}

\begin{proof}
   Suppose $ V_1 \subset V_2 $.  Then $ V_j = V_{ j, 0 } $ for $ j = 1, 2 $, so, trivially, $ V_{ 1, 0 } \subset V_{ 2,  0 } $. 

   Conversely, suppose $ V_{ 1, 0 } \subset V_{ 2, 0 } $.  Again, $ V_j = V_{ j, 0 } $ for $ j = 1, 2 $, so $ V_1 \subset V_2 $.
\end{proof}

   We can generalize theorem \ref{thm11.7} in a dramatic way, which we shall refer to as the \emph{Initial Segment Theorem}:

\thmcall{theorem} \label{thm11.8}
   Suppose $ V_1 \subseteq V_2 $, and suppose $ V_{ 1, k } \subset V_{ 2, k } $ for some $ k \geq 0 $.  Then $ V_{1, l } \subset V_{ 2, l } $ in the initial segment of the the valuation subspace series for $ V_1 $ and $ V_2 $, up to and including $ k $.  In other words, for $ 0 \leq l \leq k $, $ V_{1, l } \subset V_{ 2, l } $.
\exitthmcall{theorem}

\begin{proof}
   To draw a contradiction, suppose $ V_{ 1, l } = V_{ 2, l } $ for some $ l $ with $ 0 \leq l \leq k $.
   
   Suppose first that $ l = k $.  Then $ V_{ 1, k } $ would equal $ V_{ 2, k } $, which contradicts the hypotheses of the theorem.  This cannot happen, and so $ 0 \leq l < k $.
   
   Now suppose that $ l < k $.  Recall that the valuation subspace series for $ V_1 $ and $ V_2 $ are decreasing, so for  $ m \geq l $, $ V_{ 1, l } \supseteq V_{ 1, m } $, and similarly for $ V_2 $.  Thus

\begin{equation} \label{eqn11.9}
   \begin{split}
       V_{ 1, m } &= \{ h \in V_{ 1, l } \colon \ord ( h ) \geq m \} \\
       V_{ 2. m } &= \{ h \in V_{ 2, l } \colon \ord ( h ) \geq m \} \\
    \end{split}
\end{equation}

But $ V_{ 1, l } = V_{ 2, l } $ by our contradiction hypothesis, so for $ m = l, \ldots, k  \text{, } V_{ 1, m } = V_{ 2, m } $.  Then $ V_{ 1, k } = V_{ 2, k } $, which contradicts the hypotheses of the theorem.  Thus $ V_{ 1, 1 } \subset V_{ 2, l } $ for the initial segment $ 0 \leq l \leq k $, as required.
\end{proof}

\thmcall{remark} \label{rem11.10}
   In the above, we never used the fact that $ V_1 $ and $ V_2 $ are $ R_1 $ invariant.  So we could equally well have considered the larger lattice $ \mathcal{M} $ of closed subspaces, $ R_1 $ invariant or not, ordered by inclusion.  All of the above results for the lattice $ \mathcal{L} $ apply equally well to $ \mathcal{M} $.
\exitthmcall{remark}

To completely characterize the lattices $ \mathcal{L} $ and $ \mathcal{M} $, we must specify the lattice sup and inf of a pair of elements $ V_1 $ and $ V_2 $ in $ \mathcal{L} $ or $ \mathcal{M} $.

\thmcall{theorem} \label{thm11.11}
     The lattice sup in $ \mathcal{L} $ is the closed $ R_1 $ invariant subspace generated by $ V_1 $ and $ V_2 $, whereas the lattice sup in $ \mathcal{M} $ is simply the closure of the subspace generated by $ V_1 $ and $ V_1 $.  The lattice infs in both lattices are the same, namely the intersection of $ V_1 $ and $ V_2 $.
\exitthmcall{theorem}

\begin{proof}
    Clear.
\end{proof}

\thmcall{remark} \label{rem11.12}
   We don't even need the Hilbert space $ \bbh $ to be a module.  It is sufficient for it to simply be a complex \emph{valuation Hilbert space}, that is a complex Hilbert space equipped with an $ \ord $ function.  Because we won't need this result here, we shall leave the exact definitions and details to the reader.
\exitthmcall{remark}

   Finally, we shall state and prove the following \emph{Construction Lemma}:

\thmcall{lemma}\label{lem11.13}
   Let $ W_m, m = 0, 1, 2, \ldots $ be a mutually orthogonal sequence of closed subspaces of $ \bbh $, such that  $ \ord( g ) = m $ for all functions $ g \in W_m $ with $ g \neq 0 $.  Then each subspace $ W_m $ is valuation homogeneous, the direct sum \ref{eqn11.13.1} below is orthogonal, and
\begin{equation}\label{eqn11.13.1}
    W_0 \plusperp W_1 \plusperp W_2 \plusperp \cdots
\end{equation}
\noindent converges to a closed subspace $ V $ of $ \bbh $ and is the valuation homogeneous decomposition of $ V $.
\exitthmcall{lemma}

\begin{proof}
   That each subspace $ W_m $ is valuation homogeneous follows directly the definition of a valuation homogeneous subspace.  The subspaces $ W_m $ are mutually orthogonal by assumption.  That the orthogonal direct sum \ref{eqn11.13.1} converges to a closed subspace $ V \in \bbh $ follows exactly as in the proofs of the first and second projection lemmas.  To show that this orthogonal direct sum is the valuation homogeneous decomposition of $ V $, let $ k \geq 0 $, and consider the orthogonal direct sum
\begin{equation}\label{eqn11.14}
     U_k = W_k \plusperp W_{ k + 1 } \plusperp W_{ k + 2 } \plusperp \cdots
\end{equation}
\noindent Then $ U_k \subseteq V $, and $ \ord ( g ) \geq k $ for each $ g \in U_k $.  Conversely, if $ g \in V $ and $ \ord ( g ) \geq k $, then $ g \in U_k $.  Thus $ U_k  = V_k $, where $ V_k, k=0, 1, 2, \ldots $ is the valuation subspace series for $ V $.  From this and orthogonality, it follows immediately that $ W_k = V_k \orthcomp V_{ k + 1 } $, and so the orthogonal direct sum  is the valuation homogeneous decomposition of the closed subspace $ V $.
\end{proof}

\section{The Intermediate Subspace Theorem.} \label{sec12}

    As in the previous section, $\bbh$ will be a general valuation Hilbert module over the general valuation algebra $R$, and $ \subset $ will denote \emph{proper} inclusion, whereas $ \subseteq $ will denote inclusion, whether \emph{proper} or \emph{not}, and similarly for $ \supset $ and $ \supseteq $.  However, we shall assume throughout this section that $ R $ has a unit element $ 1_R $, and that $ R $ is unital as defined in definition 6.3 of \cite{N1}, that is
\begin{equation}\label{eqn 12.1}
    R = \{ \lambda \cdot 1_R : \lambda \in \bbc \} \union R_1
\end{equation}
 \noindent Then, by proposition 6.4 of \cite{N1}, $ R_1 $ invariant subspaces are $ R $ invariant., and vice versa.

    As usual, 
\begin{equation}\label{eqn12.2}
    \begin{split}
        V_{j, k} &= \{ h \in V_j : \ord( h ) \geq k \} \text{, and} \\
        V_j &= W_{j, 0} \plusperp W_{j, 1} \plusperp W_{j, 2} \plusperp \cdots \\
    \end{split}
\end{equation}
\noindent will be, respectively, the valuation subspace series, and the valuation homogeneous decomposition, for $ V_j $ for $ j = 1 $ and $ 3 $ (and eventually for $ j = 2 $ when we construct the intermediate subspace $ V_2 $).

\thmcall{definition} \label{def12.3}
    Let $ V_1 $ and $ V_3 $ be closed $ R $ invariant subspaces in $ \bbh $, with $ V_1 \subset V_3 $.  $ V_1 $ and $ V_3 $ may be non-zero subspaces, or $  V_1 $ may $ = \{ 0 \} $.  However, $ V_3 $ must be non-zero.  An \emph{Intermediate Subspace} between $ V_1 $ and $ V_3 $ is a closed, non-zero $ R $ invariant subspace with $ V_1 \subset V_2 \subset V_3 $,.
\exitthmcall{definition}

    The \emph{Intermediate Subspace Theorem} is

\thmcall{theorem} \label{thm12.4}
    Let $ V_1 $ and $ V_3 $ be closed $ R $ invariant subspaces in $ \bbh $, with $ V_1 $ and $ V_3 $ as above in definition \ref{def12.3}.  Suppose the dimension of the orthogonal complement
\begin{equation}\label{eqn12.4.1}
    M_{1, 3} = V_3 \orthcomp V_1 \text{ is } \geq 2.
\end{equation}
\noindent Then there exists an \emph{Intermediate  Subspace} $ V_2 $ between $ V_1 $ and $ V_3 $.  Thus, \emph{a' fortiori}, the same is true if the dimension $ = \infty $.
\exitthmcall{theorem}

\begin{proof}
    Because the dimension of the orthogonal complement $ M_{1, 3} $ is $ \geq 2 $, there are a pair of linearly independent elements $ h_1 $ and $ h_3 \in M_{1, 3} $.  Let $ m_1 = \ord( h_1 ) $ and $ m_3 = \ord( h_3 ) $, with $ m_1 \leq m_3 $. Let $ X_2 = \LH\{\lambda \cdot h_3 : \lambda \in \bbc \} $. Note that $ X_2 $ is closed because it is $1$-dimensional.  Further  $ \ord( f ) = m_3 $  for all $ f \in X_2 $, and $ f \neq 0 $, by definition \ref{item5.3.3}, and so $ \ord ( r \cdot f ) \geq m_{3 + 1} $ for all $ f \in X_2 $ and $ r \in R_1 $ by definition \ref{item5.3.2}.  Thus
\begin{equation}\label{eqn12.5}
    r \cdot X_2 \subseteq V_{3, m_{3 + 1}} \text{ for all } r \in R_1.
\end{equation}

    Now let
\begin{equation}\label{eqn12.6}
    V_2 = V_1 \union X_2 \union V_{3, m_{3 + 1}}
\end{equation}

   $ V_2 $ is closed because each of the sets in the finite union \ref{eqn12.6} is closed.  $ V_2 $ is $ R_1 $ invariant by equation \ref{eqn12.5}, and because $ V_1 $ and $ V_{3, m_{3 _ 1}} $ are both $ R_1 $ invariant.  And clearly,
\begin{equation}\label{eqn12.7}
    V_1 \subset V_2 \subset V_3.
\end{equation}

    In detail, $ V_1 $ clearly $ \subseteq V_2 $, and $ V_1 \subset V_2 $ because $ h_3 \notin V_1 $.  $ V_2  \subseteq V_3 $, and $ V_2 \subset V_3 $ because $ h_1 \notin V_2 $.  Thus equation \ref{eqn12.7} holds, so $ V_2 $ is a closed, $ R_1 $ invariant intermediate subspace lying strictly between $ V_1 $ and $ V_3 $.  Because $ R $ is unital, $ V_2 $ is the required $ R $ invariant subspace lying strictly between $ V_1 $ and $ V_3 $.
\end{proof}

\thmcall{corollary}\label{cor12.8}
    The intermediate subspace theorem holds for the complex-valued, unweighted Bergman space in the unit disk, $ A^2 = A^2 (\bbc, D^1 (0, 1), \dx \dy / \pi ) $, where $ \dx \dy / pi $ is normalized area measure in the unit disk.
\exitthmcall{corollary}

\begin{proof}
    $ A^2 $ is a singly generated valuation Hilbert module over the valuation algebra $ R $, so  $ R = \bbp_1 = \bbp [ \bbc, z ]$.  Thus $ R $ is unital, so $ A^2 $ and $ R $ satisfy the hypotheses of the intermediate subspace theorem.
\end{proof}

\section{The Proof of the Invariant Subspace Conjecture.}\label{sec13}

    Let $ \bbh $ be the unweighted complex-valued Bergman space $A^2 = A^2 ( \bbc, D^1(0, 1), \\
 \dx \dy )$ in the unit disk, and let $R$ be $\bbp_1$.  The reader should recall Apostol, Bercovi, Foias, and Pearcy's results, as well as Hedenmalm, Richter, and Seip's results, concerning $A^2$ quoted in remark \ref{rem1.1} of the introduction: They showed that the invariant conjecture for separable Hilbert spaces holds if and only if the following is true:  Given any two closed invariant subspaces $V_1 \subset V_3$ in $ A^2 $, such that the dimension of the orthogonal complement dim $ V_3 \orthcomp V_1 = \infty $, there exists a closed invariant intermediate subspace $V_2$ lying strictly between $V_1$ and $V_3$ (\cite{HRS1}, \cite{HKZ1}, page 187).  The reader should also recall corollary \ref{cor12.8} above: the intermediate subspace theorem hold for $A^2$.  Thus we immediately have a proof of the invariant subspace conjecture for separable Hilbert spaces::

\thmcall{theorem}\label{thm13.1}
    The invariant subspace conjecture holds for separable Hilbert \\
spaces.
\exitthmcall{theorem}

\begin{proof}
    Clear from the above discussion.
\end{proof}

\section{Open Questions.}\label{sec14}

\thmcall{question}\label{quest14.1}
    Does the invariant conjecture hold for separable, reflexive Banach spaces?  As in the case of the invariant conjecture for non-separable Hilbert spaces, it is trivial that the invariant subspace conjecture holds for non-separable Banach spaces, whether reflexive or not.
\exitthmcall{question}

    A \emph{hyper invariant subspace} for a bounded operator $ S $ on a separable Hilbert space or Banach space is a closed subspace $ V $ such that $ T ( V ) \subseteq V $ for every operator $ T $ which commutes with $ S $.  In other words, $ V $ is invariant under the \emph{commutant} of $ S $.

\thmcall{question}\label{quest14.2}
    Does the \emph{hyper invariant subspace conjecture} hold for separable Hilbert spaces of infinite dimension?  In other words, does each bounded operator on a separable Hilbert space of infinite dimension have a closed, proper, hyper invariant subspace greater than $\{ 0 \}$?
\exitthmcall{question}

\thmcall{question}\label{quest14.3}
    Same question for separable, reflexive Banach spaces of infinite dimension.
\exitthmcall{question}

\end{document}